\documentclass[11pt]{article}
\usepackage{amsmath}
\usepackage{amssymb}
\usepackage{amsthm}
\usepackage{graphicx}
\usepackage{psfrag}
\usepackage{verbatim}
\def\E{\mathbb{E}}

\def\eps{\epsilon}

\def\1{\mathbf{1}}

\def\tce{t_c + \eps}
\def\tce2{t_c + \frac{\eps}{2}}
\def\ER{Erd\H{o}s-R\'{e}nyi }

\newtheorem{thm}{Theorem}
\newtheorem{lem}{Lemma}

\newtheorem{defn}{Definition}
\newtheorem{prop}{Proposition}

\newtheorem{question}{Question}

\begin{document}

\title{Random $k$-SAT and the Power of Two Choices}
\author{Will Perkins\footnote{School of Mathematics, Georgia Tech, 686 Cherry St, Atlanta, GA 30332, E-mail: perkins@math.gatech.edu. Supported in part by an NSF Postdoctoral Fellowship.}}

\maketitle

\begin{abstract}
We study an Achlioptas-process version of the random $k$-SAT process: a bounded number of $k$-clauses are drawn uniformly at random at each step, and exactly one added to the growing formula according to a particular rule.  We prove the existence of a rule that shifts the satisfiability threshold. This extends a well-studied area of probabilistic combinatorics (Achlioptas processes) to random CSP's.  In particular, while a rule to delay the 2-SAT threshold was known previously, this is the first proof of a rule to shift the threshold of $k$-SAT for $k \ge 3$.  

We then propose a gap decision problem based upon this semi-random model.  The aim of the problem is to investigate the hardness of the random $k$-SAT decision problem,  as opposed to the problem of finding an assignment or certificate of unsatisfiability.  Finally, we discuss connections to the study of Achlioptas random graph processes.

\end{abstract}

\section{Introduction}

The mathematical study of phase transitions and threshold behavior in random structures began with Erd\H{o}s and R\'{e}nyi's  paper on the evolution of random graphs \cite{erdős1960evolution}.  In it, they showed that for any fixed $\eps >0$, if the number of edges in the random graph is at most $(1/2 -\eps) n$ then the largest connected component is of size $O(\log n)$ with probability $1-o(1)$ (whp), and if the number of edges is at least $(1/2 + \eps) n$, then the largest component is of size $\Theta(n)$ whp.  The existence of a giant, linear-sized component thus exhibits what is known as a sharp threshold, with its probability rising from near $0$ to near $1$ with the addition of a sublinear number of additional edges.  This is as opposed to the coarse threshold for the presence of a triangle in a random graph, a property whose probability is strictly bounded away from $0$ and $1$ for any linear number of edges.  Since that initial paper, the random graph phase transition has been studied in great detail; it is now known that the scaling window has width $n^{2/3}$ \cite{bollobas1984evolution}, the structure of both the giant component and the smaller components is well-understood, and many modifications of the original model have been studied.  The \ER random graph still plays a central role in discrete probability, both as an object in its own right and as a tool to solve other problems.

In theoretical computer science, the threshold phenomenon that has attracted the most study is the unsatisfiability threshold in the random $k$-SAT model.  An instance of random $k$-SAT on $n$ variables and $m$ clauses consists of the conjunction of $m$ clauses, each of which is the disjunction of $k$ literals and chosen uniformly at random from the set of $\binom n k 2^k$ possible clauses.  Small variations in the description of the model, such as the difference between adding clauses with or without replacement are not significant with respect to the threshold behavior.  

In \cite{friedgut1999sharp} Friedgut proved that the satisfiability threshold is sharp.  In particular, he proved that there exists a sequence $r_k(n)$ so that for every $\eps >0$,
\begin{align*}
\Pr \left[ \Phi_{(r_k(n) +\eps) n}^{(k)} \text { is satisfiable} \right] &= o(1) \\
\Pr \left[ \Phi_{(r_k(n) - \eps) n}^{(k)} \text { is satisfiable} \right] &= 1 -o(1)
\end{align*}
where $\Phi_{r n}^{(k)}$ is a uniformly random $k$-SAT formula with $rn $ clauses.  It remains an open problem whether or not the sequence $r_k(n)$ has a limit; this is known as the `satisfiability threshold conjecture'.  Much work has been done on proving upper and lower bounds on $r_k(n)$: the current best upper and lower bounds for $k=3$ are 4.508 \cite{dubois1997general} and 3.52 \cite{kaporis2003selecting}, \cite{hajiaghayi2003satisfiability} respectively. See \cite{achlioptas2009random} for best current bounds for other values of $k$ and a survey of the problem.  Typically upper bounds come from a variant of the first-moment method, while lower bounds come from either analyzing an algorithm \cite{broder1993satisfiability}, \cite{chao1986probabilistic}, \cite{achlioptas2000setting} or the second-moment method \cite{achlioptas2003threshold},\cite{achlioptas2006random}.

A second central open question in this area is whether or not random $k$-SAT formulae ($k \ge 3$) at or near the satisfiability threshold are computationally hard (2-SAT is solvable in polynomial time in the worst case). Selman, Mitchell and Levesque \cite{selman1996generating} gave experimental evidence that near the threshold certain algorithms cannot determine satisfiability efficiently.  But in some sense, the decision version of the random $k$-SAT problem is trivially easy.  Above the threshold density we can say `unsatisfiable' and below the threshold say `satisfiable' and be correct with high probability without even looking at the sampled instance.  Because of this, the study of the hardness of random $k$-SAT has turned towards refuting satisfiability for clause densities above but as close as possible to the threshold (eg.  \cite{feige2004easily}), and towards finding satisfying assignments at as high densities as possible.  

In this work we consider a variant of the random $k$-SAT problem for which the decision problem is still relevant.  Our initial motivation came from the work on Achlioptas processes in the study of random graphs.  Dimitris Achlioptas initiated this study by asking whether, given the choice of two random edges at each step of a random graph process, one could delay the phase transition by a constant factor.  His question was in fact motivated by understanding algorithms for random $3$-SAT, and in particular understanding how to choose a variable to set in a free step of a non-backtracking algorithm so an underlying $2$-SAT formula remains satisfiable \cite{AchEmail}.    

Formally an Achlioptas random graph process is defined as follows:
\begin{itemize}
\item Begin at step $0$ with an empty graph on $n$ vertices.
\item At step $i$, two uniformly random edges are presented, and exactly one of the two is selected according to a given rule and added to the graph
\item The choice of edge can depend on the edges presented, the current graph, and the history of the process, but not on the edges to be presented in subsequent steps.
\end{itemize}
 His original question was whether there is a rule for choosing one of  the two edges so that at step $(1/2+ \eps) n$ the graph contains no linear-sized connected component whp.  His question was answered affirmatively by Bohman and Frieze \cite{bohman2001avoiding} in the first of many papers to study Achlioptas processes.  While the phase transition has received the most attention (\cite{bohman2006creating}, \cite{flaxman2005embracing}, \cite{riordan2012achlioptas}), rules for shifting the threshold of other properties have also been found (e.g. Hamiltonian cycles \cite{krivelevich2010hamiltonicity} or small subgraphs \cite{krivelevich2009avoiding}, \cite{mutze2011small}).  One primary aim of the study of Achlioptas processes is to understand which qualitative properties of the phase transition are robust under small modifications of the model.  In \cite{janson2010phase}, \cite{kang2012bohman}, \cite{bhamidi2013aggregation} it is shown that certain critical exponents are universal for a large class of Achlioptas rules.

Sinclair and Vilenchik \cite{sinclair2010delaying} first considered the Achlioptas process model with regard to a random CSP.  In particular, they exhibited a rule for choosing one of two random clauses that delays unsatisfiability for random 2-SAT by a constant factor.  They also considered `off-line' rules for random $k$-SAT and on-line rules for $k$-SAT with $k = \omega (\log n)$.

In this work we study the Achlioptas-process version of random $k$-SAT (see Section \ref{modelSec} for a formal definition), and show in Section \ref{ResultsSec} that in fact the satisfiability threshold can be shifted, for any $k$, and in particular for the computationally interesting cases $k \ge 3$.  Studying this semi-random model of $k$-SAT is a step towards understanding the standard model better, as has been done in the case of random graphs. 

To address the question of the hardness of the random $k$-SAT decision problem while avoiding the triviality of almost sure satisfiability below the threshold and almost sure unsatisfiability above the threshold, we use the semi-random model of $k$-SAT described above and propose an accompanying decision problem in Section \ref{GapSec}.

To sum up, this paper makes the following contributions:
\begin{enumerate}

\item We prove that a specific Achlioptas rule shifts the satisfiability threshold of random $k$-SAT by a constant factor for all $k \ge 2$.  This was previously known only for the case of 2-SAT. 
\item In particular we show that there is a 5-clause rule to shift the 3-SAT threshold.
\item We improve the previously-known factor of delay for a 2-clause, 2-SAT rule.
\item The rule we analyze and our method of proof also shows that biased random $k$-SAT formulas are easier to satisfy than unbiased ones.  
\item We propose a gap decision problem for this semi-random Achlioptas model that may be more approachable than random $k$-SAT from the perspective of computational complexity.
\end{enumerate}

\section*{Notation}

The set of binary variables on which our random formulae are built is $\{x_1, \dots x_n\}$, and all asymptotics are as $n \to \infty$.  A literal is a variable $x_i$ or its negation $\overline x_i$, and we will denote literals with the letter $w$.  A $k$-clause is a disjunction of $k$ literals, $(w_{i_1} \vee w_{i_2} \vee \cdots \vee w_{i_k})$. A formula of $m$ clauses is the conjunction of $m$ $k$-clauses and is satisfiable if there exists an assignment to the $n$ variables that satisfies each of the $m$ clauses.  We will write $\Phi_m$ for a formula of $m$ clauses.  We write that an event $E$ holds with high probability or whp if $\Pr[E] \to 1$ as $n \to \infty$.  

\section{Semi-Random $k$-SAT: The Model}
\label{modelSec}

Here we define an $l$-clause Achlioptas $k$-SAT process analogously to an Achlioptas random graph process:
\begin{enumerate}
\item Begin at step $0$ with an empty formula, $\Phi_0 = \emptyset$.
\item Each each step, $l$ clauses are selected uniformly at random, with replacement, from the $\binom n k 2^k$ possible $k$-CNF clauses.  
\item According to a fixed rule $R$, exactly one of the $l$ clauses is chosen and added to the current formula.  $\Phi_i = \Phi_{i-1} \wedge \phi_i$, where $\Phi_{i-1}$ is the current formula, and $\phi_i$ is the clause chosen at step $i$.
\end{enumerate}
Note that different rules $R$ lead to a different processes (and different distributions over formulas $\Phi_m$ at step $m$).  The rule `Always select the first clause' leads to the classic random $k$-SAT distribution.  The rule $R$ can be a function of the $l$ clauses presented, the current formula $\Phi_{i-1}$, and the entire history of presented clauses up to step $i$.  The rule can also use randomness.  The rule, however, cannot be a function of the clauses presented in subsequent steps (such rules, while not standard Achlioptas processes,  are called `off-line' rules and have been studied for both the random graphs \cite{bohman2004avoidance} and $k$-SAT formulae \cite{sinclair2010delaying}).  We do not require the rule to have an efficient implementation; a rule such as `Add the first clause unless it creates an unsatisfiable formula' is allowed, and the resulting process would resemble the random satisfiable k-CNF process from \cite{krivelevich2009random}. The rules we analyze below will be simpler: they will be functions of only the $l$ clauses presented in the current round.

\section{Results}
\label{ResultsSec}

Let $r_k=r_k(n)$ be the threshold for the satisfiability of random $k$-SAT. Our first result is that the satisfiability threshold can be delayed with an $l$-clause rule, for $l$ at most $5$.  
\begin{thm}
\label{GenDelThm}
For every integer $k \ge 2$ and some fixed $\eps >0$, there exists a $l$-clause Achlioptas rule for random $k$-SAT so that with probability $1-o(1)$, the formula generated after $(r_k + \eps) \cdot n $ steps is satisfiable.    For $k\ge 7$, we can take $l=3$; for $3\le k \le  6$, we can take $l=5$; and for $k=2$, we can take $l=1$.
\end{thm}

For the case $k=2$, we improve the constant factor of delay for a 2-clause rule in the results of \cite{sinclair2010delaying} with a different rule and a different proof.
\begin{thm}
\label{k2Thm}
There is a 2-clause Achlioptas rule for random 2-SAT that generates a formula that is satisfiable whp after $1.055 n$ steps.
\end{thm}

The proofs  follow in Section \ref{ProofSec}.

\section{Semi-Random Gap $k$-SAT}
\label{GapSec}

Once we know that an Achlioptas rule can change the satisfiability threshold, the following decision problem becomes meaningful.

Let $\Phi^{(R)}_i$, $i=1, \dots$ be a growing semi-random 3-SAT formula generated according to some 2-clause Achlioptas $3$-SAT rule $R$.  As an input, we are presented with the sequence of clauses, but not the rule $R$.  We want to distinguish between the following two cases:

\begin{itemize}
\item NO: At step $4 n$, the formula is unsatisfiable.
 \item YES: At step $5 n$, the formula is satisfiable.
\end{itemize}
We say the algorithm is in error if $\Phi^{(R)}_{4n}$ is unsatisfiable and the algorithm returns `YES' or if  $\Phi^{(R)}_{5n}$ is satisfiable and the algorithm returns `NO'.  If the first unsatisfiable step occurs between $4n$ and $5n$, then either answer is acceptable.

\begin{question}
 Is there an efficient algorithm $A$, so  that for all rules $R$, 
 \[ \Pr_A[\text{error}] = o(1)  ?\] 
\end{question}

In other words, we are asking for an algorithm that gives an acceptable answer with high probability over the randomness in the clauses presented, in the worst case over all possible rules $R$.  

We can generalize the above problem to $k-$SAT and by varying three parameters: the number of random clauses presented at each step $(l)$, the lower threshold ($c_1$), and the upper threshold ($c_2$).  We have chosen $k=3, l=2, c_1= 4, c_2=5$ for simplicity.  The problem becomes harder as $l$ increases, and easier as either $c_2 $ increases or $c_1$ decreases.  If the adversary (i.e. the rule $R$) had no choice, and a random clause was added at each step, the problem would be easy: the satisfiability threshold would occur at $r_k(n)$, and if $r_k(n) \le c_1$, NO would be acceptable whp; if $r_k(n) \ge c_2$ YES would be acceptable; and otherwise either answer would be acceptable.  If the rule could choose from all $2^k \binom n k$ possible clauses at each step, the problem would be NP-hard.  This semi-random decision problem interpolates between the two extremes.

One way to interpret Achlioptas' original question is whether, under a specific model of semi-random graphs, the phase transition occurs when the average degree of a vertex hits $1$ as it does in the \ER random graph.  Bohman and Frieze answered `no' to this question:  in a sense the location of the phase transition at  average degree $1$ in the standard model is an artifact of the independence and uniformity of the random edges.  However, the study of Achlioptas processes has identified a different statistic, rather than the average degree, that does control the phase transition for a large class of  processes.  This is the \emph{susceptibility}, or the average component size in the graph: $S(G) = n^{-1} \sum_{v} |C(v)|$.  Bohman and Kravitz \cite{bohman2006creating} and Spencer and Wormald \cite{spencer2007birth} show that for the class of `bounded-size' rules, the blow-up point of an ODE tracking the growth of $S(G)$ marks the critical point for the phase transition.  

The susceptibility allows one to understand where and why the phase transition occurs in Achlioptas random graph processes but is not needed algorithmically, since detecting a giant component is an easy computational problem.  In the case of $k$-SAT however, we can ask if there is such a statistic, correlated with the threshold, which is efficiently computable.

\section{Proofs}
\label{ProofSec}

Theorems \ref{GenDelThm} and \ref{k2Thm} are corollaries of the following lemma:

\begin{lem}
\label{MainLem}
 For fixed integers $k \ge 2$ and $l \ge 2$, there exists an $l$-clause Achlioptas rule for random $k$-SAT which creates a formula that, for every $\eps >0$, is satisfiable whp after $(r(k,l) -\eps) n$ steps, where 
\begin{equation}
\label{rklEQ}
 r(k,l) =  \frac{1}{\left( \frac{k+1}{2^k}  \right)^ {l-1} \frac{k}{2^k}   + 2 \sqrt{ \left( \frac{k+1}{2^k}  \right)^ {l-1}  \frac{1}{2^k} \left( 1- \left( \frac{k+1}{2^k}  \right)^ l \right)   }   }
\end{equation}
\end{lem}

To prove Theorem \ref{GenDelThm}, we consider a few cases.  For $k=3$, we have $r(3,5) = 5.06508\dots$, which is strictly larger than the best current upper bound on the random $3$-SAT threshold, so $l=5$ choices is enough.  For $k=4,5,6$ we check numerically that $r(k,5)\ge 2^k \ln 2$, an upper bound on the satisfiability threshold.  For $k=7$, we check numerically that $r(k,3) > 2^k \ln 2$, and the function
\[ g(k) := \frac{ 2^k \ln 2}{ r(k,3)} \]
is decreasing in $k$, so $3$ choices suffice for all $k \ge 7$.  

For Theorem \ref{k2Thm}, we plug in $k=2$ and $l=2$, and get $r(k,l) = 1.05505\dots > 1.055$.  


To prove Lemma \ref{MainLem}, we start by giving an idea of the strategy we will use to delay unsatisfiability.  Consider a biased random 3-SAT formula generated by adding $m= rn$ clauses, with each clause selected as follows: with probability $p$, choose a clause uniformly from all 3-clauses with $3$ positive literals, and with probability $1-p$ choose a clause uniformly from all 3-clauses.  Let $\Phi_{rn}$ denote that random formula.  This formula is biased in favor of assignments with more $+1$'s than $-1$'s.  Let $Z_\beta$ be the number of satisfying assignments to $\Phi_{rn}$ with $\beta n$ $+1's$ and $(1-\beta) n$ $-1$'s, and let $x_\beta$ be the particular assignment that assigns the first $\beta n$ variables $+1$ and the rest $-1$.  Then
\begin{align*}
\frac{1}{n} \log  \E Z_\beta &=\frac{1}{n} \log \left(  \binom{n}{\beta n} \Pr[ x_\beta \text { satisfying} ]  \right )\\
&= H(\beta) + r \log \Pr[ x_\beta \text{ satisfies }\Phi_1  ] + o(1)
\end{align*}
where $H(\beta)$ is the binary entropy function.
\begin{equation*}
\Pr[ x_\beta \text{ satisfies }\Phi_1  ] \sim p (1- (1-\beta) ^3) + (1-p) \cdot 7/8
\end{equation*}
If we pick $r$ small enough so that $\max_{\beta \in (0,1)} H(\beta) + r  \log  [ p (1- (1-\beta) ^3) + (1-p) \cdot 7/8] >1$, then there will be exponentially many satisfying assignments in expectation.  One can show that for any $p >0$, there is some such $r $ larger than $\log 2 / \log (7/8) \sim 5.19...$ which is the simple first-moment upper bound for random 3-SAT.  And in fact for any $r >0$, there is a $p \in (0,1)$ so that $\max_{\beta \in (0,1)} H(\beta) + r \log [ p (1- (1-\beta) ^3) + (1-p) \cdot 7/8] >1$, which shows that with enough bias, the first-moment bound can be pushed arbitrarily high.  This complements results on random regular $k$-SAT \cite{boufkhad2006regular}, in which an extreme lack of bias, with every literal having the same degree, leads to an earlier threshold.  

Having exponentially many solutions in expectation does not imply a single solution with significant probability, so to prove the lemma we will use a related but different rule. This rule will select one of the $l$ clauses presented at each step as follows:
\begin{itemize}
 \item If one of the first $l-1$ clauses contains at least two positive literals, add it. (If there is more than one such clause, add the first.)
\item Otherwise add clause $l$.
\end{itemize}

The effect of this rule is to bias the formula in favor of majority $+1$ assignments as above. Let $r = r(k,l) -\eps$. To prove that this rule produces a satisfiable formula whp, we will begin by taking the formula at step $r n$ and converting it into a 2-SAT formula.  For each clause with two or more positive literals, we keep a 2-CNF clause with only the first two positive literals; for each clause with exactly one positive literal, we keep a 2-clause with the positive literal and the first negative literal; and for each clause with all negative literals, we keep a 2-clause with the first two literals.  This gives a 2-SAT formula with $r n$ clauses, call it $F_2(r n)$.  If $F_2(rn)$ is satisfiable, the original formula is also satisfiable since each 2-clause is a sub-clause of the corresponding $k$-clause.  Each clause is also added independently, and the probability it has two, one or zero positive literals is, respectively:
\begin{align*}
p_0 &= \Pr \left [Bin \left(k,1/2 \right) \le 1 \right]^{l-1} \cdot\Pr \left [Bin \left(k,1/2 \right) =0 \right]  \\
&= \left( \frac{k+1}{2^k}  \right)^ {l-1} \cdot \frac{1}{2^k}\\
p_1 &=  \Pr \left [Bin \left(k,1/2 \right) \le 1 \right]^{l-1} \cdot\Pr \left [Bin \left(k,1/2 \right) =1 \right]  \\
&= \left( \frac{k+1}{2^k}  \right)^ {l-1} \cdot \frac{k}{2^k} \\
 p_2 &=  1-\Pr \left [Bin \left(k,1/2 \right) \le 1 \right]^{l}  \\
 &=1- \left( \frac{k+1}{2^k}  \right)^ l
\end{align*}
Also, each clause with $i$ positive literals is distributed uniformly from all 2-clauses with $i$ positive literals. As is standard in studying random graphs and random CSP's, we will consider the model of random 2-CNF formulas in which each of the $\binom n 2$ possible clauses with two positive literals is present in a random formula $\overline F_2(r,n)$ independently with probability $q_2$, each of the $n(n-1)$ clauses with one positive and one negative literal is present with probability $q_1$, and each of the $\binom n 2$ with two negative literals is present independently with probability $q_0$.  If we set
\begin{align*}
q_2 &= \frac{ 2 p_2 r}{n}  \\
q_1 &= \frac{p_1 r }{n}   \\
q_0 &=  \frac{2 p_0 r}{n}  
\end{align*}
then proving satisfiability whp of $\overline F_2(r,n)$ implies satisfiability whp of $F_2(rn)$ (see \cite{bollobas2001scaling}, Appendix A, for details of the equivalent behavior of the two models). 

To study the satisfiability of $\overline F_2(r,n)$ we will use an approach of \cite{chvatal1992mick} and \cite{cooper2002note}, following \cite{aspvall1979linear}, where it is noted that if there is no bicycle in a formula's `implication graph', the formula is satisfiable.  The vertices of the implication graph are the $2n$ literals, and for each clause $(w_i \vee w_j)$ in the formula, we add two directed edges, $(\overline w_i \to w_j)$ and $(\overline w_j \to w_i)$ to the graph. 
\begin{defn}
A bicycle of length $t$ is a sequence of $t$ literals of distinct variables, $w_1, w_2, \dots w_t$ where the $t-1$ directed edges $(w_1 \to w_2), (w_2 \to w_3), \dots (w_{t-1} \to w_t)$ are present in the implication graph, as well as two additional directed edges $(u \to w_1)$ and $(w_t \to v)$ for $u,v\in \{w_1, \dots w_t, \overline w_1, \dots \overline w_t \}$. 
\end{defn}

\begin{prop}\cite{chvatal1992mick}
If a 2-SAT formula contains no bicycle, then it is satisfiable.
\end{prop}

To show that there are no bicycles, we proceed as in \cite{cooper2002note}, and consider first directed paths of length $\ge  L = K \eps^{-1} \log n$.  Note that a clause with two positive literals adds two directed edges from negative literals to positive literals; a clause with one negative and one positive literal adds one directed edge from a positive to a positive literal and one from a negative to a negative; and a clause with two negative literals adds two directed edge from positive literals to negative literals.  So in a path of $L$ literals, if the signs of the literals switch along the path $i$ times there must be $L-1 -i$ corresponding clauses with exactly one positive literal, and  $ i$ corresponding clauses with two or zero positive literals each, at least $\lfloor i/2 \rfloor \ge (i-1)/2$ of which must have zero positive literals.  The probability that there is a directed path connecting a given sequence of $L$ literals with $i$ sign changes is therefore bounded above by $q_1^{L-1-i} q_2^{(i+1)/2} q_0^{(i-1)/2}$, since $q_0 \le q_2$. If we let $Y_L$ denote the number of directed paths of length $L$, then, summing over the number of sign switches,
\begin{align}
\E Y_L &\le  \sum_{i=0}^L 2 n^L \binom{L-1}{i} q_1^{L-1-i} q_2^{(i+1)/2} q_0^{(i-1)/2}  \\
&= 2 n^L \sum_{i=0}^L \binom {L-1} {i} \frac{r^{L-1}}{n^{L-1}} p_1^{L-1-i} (2p_2)^{(i+1)/2} (2p_0)^{(i-1)/2}\\
&= 2 n r^{L-1} \sqrt{\frac{p_2}{p_0}} \sum_{i=0}^L \binom {L-1} {i}  p_1^{L-1-i} (2 \sqrt{p_0 p_2})^i\\
\label{pathsINEQ}
&= 2 n  \sqrt{\frac{p_2}{p_0}} r^{L-1} (p_1 + 2\sqrt{p_0 p_2})^{L-1}
\end{align}

Now let $Z_L$ be the number of bicycles of length at most $L$. Then,
\begin{align}
\E Z_L &\le \sum_{t=2}^L 2n^t (q_2 + q_1+q_0)^2 t^2 \sum_{i=0}^{t-1} \binom {t-1} i q_1^{t-1-i} q_2^{(i+1)/2} q_0^{(i-1)/2}\\
&= \sum_{t=2}^L 2n^t (2p_2 + p_1+2p_0)^2 t^2  \frac{r^{t+1}}{n^{t+1}} \sum_{i=0}^{t-1} \binom {t-1} i p_1^{t-1-i} (2p_2)^{(i+1)/2} (2p_0)^{(i-1)/2}\\
&= \sum_{t=2}^L 2n^t (2p_2 + p_1+2p_0)^2 t^2  \sqrt{\frac{p_2}{p_0}}  \frac{r^{t+1}}{n^{t+1}} \sum_{i=0}^{t-1} \binom {t-1} i p_1^{t-1-i} (2 \sqrt{p_0 p_2})^i\\
\label{BicINEQ}
&\le\frac{8}{n}\sqrt{\frac{p_2}{p_0}}  \sum_{t=2}^L   t^2 r^{t+1} (p_1 + 2\sqrt{p_0 p_2})^{t-1}
\end{align}

Considering (\ref{pathsINEQ}) and (\ref{BicINEQ}), if 
\begin{equation}
\label{RINEQ}
r \le \frac{1}{p_1 + 2\sqrt{p_0 p_2}} - \eps
\end{equation}
 for some fixed $\eps >0$, then for $L = \Theta_{\eps}(\log n)$, we have $\E Y_L = o(1)$ and $\E Z_L =o(1)$, and therefore the formula is satisfiable whp. Substituting the values of $p_0, p_1, p_2$ into (\ref{RINEQ}) yields (\ref{rklEQ}). 


\section{Discussion and Open Problems}

We conclude with some remarks and open problems.

A first question is whether there is a 2-clause rule to shift the $k$-SAT threshold for $k>2$.  It is natural to conjecture that in fact the rule used in the proof of Theorem \ref{GenDelThm} does in fact shift the threshold for $l=2$, since it does shift the first-moment upper bound.  The difficulty lies in the gap between the current upper and lower bounds on $r_k(n)$ --- to prove that the rule has shifted the threshold for 3-SAT, we need to prove that it has shifted the threshold all the way past $4.508$.  What would be much more straightforward to prove would be that certain algorithms (the unit-clause algorithm and its variants described in \cite{achlioptas2001lower}) succeed at higher densities for this rule than for random $k$-SAT.

\textit{Remark:}  Building on the ideas in this work Dani, Diaz, Hayes and Moore \cite{DaniDHM13} studied the problem and by refining the strategy and analysis, lowered the number of choices needed to delay the satisfiability threshold, in particular showing that $2$ choices suffice for $3 \le k \le 25$.   

We have discussed bounds on the satisfiability threshold of Achlioptas processes here and mentioned Friedgut's result on the sharpness of the $k$-SAT threshold.  In fact the rule we analyze can be shown to have a sharp threshold using Bourgain's sharp threshold criterion (Bourgain's appendix to \cite{friedgut1999sharp}). It would be interesting to determine which Achlioptas $k$-SAT processes have a sharp threshold or if all rules for a fixed $l$ have a sharp threshold.

\begin{question}
For fixed $l$ and $k$, is there an $l$-clause Achlioptas rule for $k$-SAT that does not have a sharp threshold?
\end{question}

 Next we note that the rules we analyze above all operate by biasing the formula to favor a particular assignment.  
 
 \begin{question}
 Can the $k$-SAT threshold be shifted by an Achlioptas rule that is symmetric with respect to assignments?
 \end{question}
 One candidate for such a rule would be the following:
 \begin{itemize}
\item If all (or none, for the opposite effect) of the literals in the first clause appear in the current formula, add it.  
\item Otherwise add the second clause.
\end{itemize}

 And finally, analogies to work done on the phase transition of Achlioptas random graphs suggest many avenues for future work.

\subsection*{Acknowledgements}
The author would like to thank Dimitris Achlioptas for explaining the motivation for his original question on the power of two choices in a random graph, and the anonymous referees for useful and constructive comments.

\bibliographystyle{abbrv}
\bibliography{sat}

\end{document}